\documentclass[12pt]{article}

\usepackage{amsmath, amssymb, amsthm}
\usepackage{xcolor}

\newtheorem{thm}{Theorem}[section]

\newtheorem{lem}[thm]{Lemma}
\theoremstyle{definition}

\newtheorem{example}[thm]{Example}
\newtheorem*{rem*}{Remark}
\newtheorem{prop}[thm]{Proposition}
\newtheorem{cl}{Claim}

\numberwithin{equation}{section}

\def\K{\mathcal{K}}
\def\R{\mathbb{R}}

\def\dom{\mathrm{dom}\,}

\def\e{\varepsilon}
\def\eps{\varepsilon}
\def\ffi{\varphi}

\def\ds{\displaystyle}

\def\ba{\begin{array}}
\def\ea{\end{array}}
\def\be{\begin{equation}}
\def\ee{\end{equation}}

\def\ra{\Rightarrow}
\def\ex{\exists}
\def\fa{\forall}
\def\nequiv{\not\equiv}

\title{Hadamard inverse mapping theorem for certain continuously G\^ateaux differentiable mappings}
\author{M. Hamamdjiev\thanks{Sofia University, Faculty of Mathematics and Informatics, \\
   5, James Bourchier Blvd, 1164 Sofia, Bulgaria, e-mail: mihailh@fmi.uni-sofia.bg}, N. Zlateva\thanks{Sofia University, Faculty of Mathematics and Informatics, \\
   5, James Bourchier Blvd, 1164 Sofia, Bulgaria, e-mail: zlateva@fmi.uni-sofia.bg}}
\date{ }

\begin{document}

\maketitle

\begin{abstract}
Hadamard  inverse mapping theorem concerns global inversion of continuously  differentiable mappings between Banach spaces. Here we consider a weaker notion of smoothness and work with  continuously G\^ateaux differentiable mappings. We  prove that a   continuously G\^ateaux differentiable mapping with invertible derivatives such that the inverses are uniformly bounded,  possesses a continuous right inverse on every compact set. We also  provide  sufficient conditions for global invertibility of certain continuously G\^ateaux differentiable mappings between Banach spaces.
\end{abstract}

\textbf{Key words: } Inverse Mapping  Theorem, continuously G\^ateaux differentiable mappings

\textbf{Mathematics Subject Classification: } 47J07, 47J30

\section{Introduction}
Invertibility of mappings is a cornerstone in  Mathematical Analysis. Let $X$ and $Y$ be Banach spaces and let
  $f : X\to Y$ be a $C^1$-smooth mapping   whose derivative $f'(x)$ is an isomorphism for every $x\in X$. Then the classical Inverse Mapping Theorem implies that $f$ is locally invertible around each point. If, in addition, $f$ satisfies the  condition
\[
\|[f'(x)]^{-1}\|\le M,\quad \forall x\in X,
\]
then $f$ is globally $C^1$ invertible on $X$, and hence a global diffeomorphism of $X$ onto $Y$. In other words, there is a $C^1$ function $g$ such that
\[
g(f(x))=f(g(x))=x,\quad \forall x\in X.
\]
This   condition    is due to Hadamard~\cite{Hadamard} for mappings between finite-dimensional spaces, and to L\'evy~\cite{Levy} for mappings between Hilbert spaces and it is known as the Hadamard-L\'evy condition.

Detailed information about this and other conditions for global inversion of smooth mappings between Banach spaces can be found, for instance, in~\cite{Plastock,Gutu}. The result  has been extended in several directions over the years. Some authors weaken the $C^1$ assumption  on the mapping using a pseudo-Jacobian, Clarke subdifferential concept instead, see~\cite{Pourciau1,Pourciau2,GutuJaramillo}; others consider set-valued maps~\cite{Pales,AubinFrankowska} or more general  underlying spaces~\cite{Ekeland,CFR}; a local surjection theorem with continuous inverse is obtained in \cite{EkelandSere}, etc.

We are interested in the global inversion of mappings between Banach spaces under smoothness assumptions weaker than $C^1$.
Mappings called continuously G\^ateaux differentiable, written  $C^1_G$, seem to be good candidates. In contrast to the $C^1$ case, a $C^1_G$ mapping may fail to be invertible on any neighborhood of a point at which its derivative  is invertible, see Example~\ref{Example}. Hence, even local invertibility of a $C^1_G$ mapping $f:X\to Y$ around a fixed point $x$ needs an additional requirement. One such   condition is
\begin{equation*}
\limsup_{x'\to x} \| f'(x') - f'(x) \|  < \|[f'(x)]^{-1} \|^{-1}.
\end{equation*}
It appears in \cite[Theorem 10]{IJLZ} and shows that a certain controlled discontinuity of the G\^ateaux derivatives  allows local invertibility.

The paper is organized as follows. Section~\ref{sec:prelim} contains the necessary definitions and preliminaries. In Section~\ref{sec:Surj} we prove that a mapping $f:X\to Y$ which is $C^1_G$ on $X$, has invertible derivatives and satisfies the Hadamard-L\'evy condition is surjective.  Section~\ref{sec:LocInv} is devoted to local invertibility of $C^1_G$ mappings. In Section~\ref{sec:Main} we prove our main result on the global  invertibility of $C^1_G$ mappings.

\section{Preliminaries and notation}\label{sec:prelim}

Let  $(X,\|\cdot\|_X )$  and $(Y,\|\cdot\|_Y )$ be Banach spaces. In $X$ we denote the unit sphere by $S_X$, and the  closed (resp. open) unit ball by $B_X $ (resp.   $B_X ^\circ $).   For a non-empty set $A\subseteq X$ we denote by $C(A,Y)$ the space of all norm-to-norm continuous mappings $f:A\to Y$. It is well-known that, equipped with the supremum norm
\[
\|f\|_A :=\sup\{\|f(x)\|_Y :x\in A\},
\]
the space $(C(A,Y),\|\cdot \|_A )$ is a Banach space provided that $A$ is compact.

A mapping $f:X\to Y$ has a  \emph{directional derivative } at a point $x\in X$ in the direction $h\in X$ whenever the limit
\be\label{dir_der}
f'(x,h):=\lim_{t\downarrow0}\frac{f(x+th)-f(x)}{t}
\ee
exists. The mapping $f$ is said to be \emph{G\^ateaux differentiable at} $x$ whenever the directional derivatives $f'(x,h)$ exist for all $h\in X$ and are bounded linear operators, denoted by $f'(x,h)=f'(x)h$. The norm of the linear map $f'(x):X\to Y$ is the standard linear operator norm    $\|f'(x)\|:=\sup \{\|f'(x)h\|_Y :\|h\|_X\le 1\}$.

A mapping $f:X\to Y$ is called \emph{Fr\'echet differentiable at}  $x\in X$ if there is a bounded linear operator $f'(x):X\to Y$ such that
\[
 \lim_{\|h\|_X\to 0}\frac{f(x+h)-f(x)-f'(x)h}{\|h\|_X}=0,
\]
and it is called \emph{smooth}, denoted $f\in C^1$, if the mapping $x\to f'(x)$ is norm-to-norm continuous.

 A mapping  $f:U\to Y$ is said to be \emph{continuously G\^ateaux differentiable} on an open set $U\subseteq X$, or $f$ is $C^1_G$ on $U$, if $f$ is continuous on $U$, the directional derivatives $f'(x,h)$ exist for all $x\in U$ and $h\in X$ and $f'(\cdot,\cdot):U\times X\to Y$ is  a \emph{jointly continuous mapping}   on the product $U\times X$. As  noted at \cite[p.73]{Hamilton} this notion does not agree with the usual definition of a $C^1$-smooth mapping in the Fr\'echet sense. If $f$ is $ C^1_G$ on $U$, then the derivative will of course be the same linear map, but the continuity requirement is weaker. It is well known, see e.g. \cite[Lemma~3.2.3, Theorem~3.2.5]{Hamilton},  that if $f$ is $ C^1_G$ on $U$, then for any $x\in U$, $f$   is G\^ateaux differentiable at $x$. It is worth  mentioning that for any $x_0\in U$ there exist $\delta >0$ and $L>0$ such that $\|f'(x)(h)\|_Y\le L\| h\|_X$ for all $x\in \{x_0\}+\delta B^\circ_X$ and $h\in X$, so $f$ is locally Lipschitz on $U$ \cite[Theorem~5.1.1]{Hamilton}.

If $f$ is $ C^1_G$ on $X$, then  $\fa x\in X$, $\fa h\in X$,
\begin{equation}\label{int2}
\frac{f(x+th)-f(x)}{t}=\int_0^1f'(x+uth)h\,du,
\end{equation}
see \cite[Lemma 3.2.4]{Hamilton}.

\bigskip

In the sequel we will make essential use of the following auxiliary result.

\begin{lem}
\label{lem:Approx}
Let  $f:X\to Y$ be   $ C^1_G $ on $X$. Let $K\subseteq X$ be a non-empty compact set, $E\subseteq X$ be a finite dimensional subspace and let $r>0$. Then there exists a function $\alpha :[-1,1]\to \R^+$ such that $\alpha(t)\to0$ as $|t|\to0$, and
\[
\sup\{\|f(x+th)-f(x)-tf'(x)h\|_Y :x\in K,\ h\in rB_E\}\le|t|\alpha(t).
\]
\end{lem}
\begin{proof}

Let $t\in \R$ be such that  $|t|\le 1$. Using \eqref{int2}, we have
\begin{eqnarray*}
\|f(x+th)-f(x)-tf'(x)h\|_Y&=&\left\|t\int_0^1 \left[f'(x+uth)h-f'(x)h\right]\, du\right\|_Y\\
  &\le& |t|\int_0^1\|f'(x+uth)h-f'(x)h\|_Y\, du.
\end{eqnarray*}
Set
\[
\alpha(t):=\sup\left\{\int_0^1\|f'(x+uth)h-f'(x)h\|_Y\,du:(x,h)\in K\times r{B}_E\right\}.
\]
Fix $\e>0$. By uniform continuity of the continuous mapping $f':X\times X\to Y$ on the compact $(K+rB_E)\times rB_E$, there is $\delta>0$ such that $\fa x_1,x_2\in K+rB_E$ and $\fa h_1,h_2\in rB_E$ satisfying  $\|x_1-x_2\|_X+\|h_1-h_2\|_X<\delta$, it holds that  $\|f'(x_1)h_1-f'(x_2)h_2\|_Y<\e$.
For $|u|\le 1$,  $|t|<1$ and $(x,h)\in K\times rB_E$ it follows that $(x+uth,h)\in (K+rB_E) \times rB_E$. Then, $\fa |t|<\min\{1,\delta/r\}$ we get
\[
\|f'(x+uth)h-f'(x)h\|_Y<\e,
\]
which implies $\alpha(t)<\e$. Therefore,  $\alpha(t)\to0$ as $|t|\to 0$.
\end{proof}

The next lemma is a precursor to the Ekeland Variational Principle, see~\cite[Chapter 5, Section~1]{Aubin-Ekeland-book}. See also the comments concerning the Basic Lemma\ on \cite[p. 93]{Ioffe-book}. It easily follows from the Ekeland Variational Principle itself, see e.g. \cite[Basic~Lemma]{Ioffe-book}, and is a partial case of Takahashi's Nonconvex Minimization Theorem~\cite{Takahashi}.   Here we state it in the form given in \cite{IZ-Hadamard}.

\begin{lem}\label{lem:mu}\cite[Lemma 3]{IZ-Hadamard}
Let $(X,\|\cdot \|_X)$ be a Banach space and $\mu:X\to\R^+\cup\{\infty\}$ be a lower semicontinuous proper function such that for some $r>0$
\[
\fa x: 0<\mu(x)<\infty \ra \ex y: \ \mu(y)<\mu(x)-r\|y-x\|_X.
\]
Then, for each $ x\in \dom\mu$ there is ${y}\in X$ such that
\[
\mu({y})=0 \text{ and } r\|{y}-x\|_X\le \mu(x).
\]
\end{lem}

\section{Surjectivity}\label{sec:Surj}

We will first prove that any $f:X\to Y$ which is $C^1_G$ on $X$ and possesses   invertible G\^ateaux derivatives satisfying the Hadamard-L\'evy condition has a continuous right inverse  on every compact set in $Y$. In Theorem \ref{thm:surj} we will use this result  to prove the surjectivity of $f$.

\begin{prop}
\label{prop:RightInv}
Let $(X,\|\cdot\|_X)$ and $(Y,\|\cdot\|_Y)$ be Banach spaces. Let  the mapping $f:X\to Y$  be  $ C^1_G$ on $X$ and possesses  an invertible G\^ateaux derivative at every point, that is $\fa x\in X$ there exists $\left[f'(x)\right]^{-1}:Y\to X$.
Let the Hadamard-L\'evy condition
\be\label{Had}
\|[f'(x)]^{-1}\|\le M,\quad \forall x\in X
\ee
hold.

Then, for every non-empty compact $\K\subseteq Y$ one can find a continuous mapping $ g\in C(Y,X)$ which is a right inverse of $f$ on $\K$, that is
$$f(g(y))=y,\quad \forall y\in \K.$$
Furthermore, if $0\in\K$ and $f(0)=0$, then $g(0)=0$.
\end{prop}

\begin{proof}
Let us fix a non-empty compact set $\K\subseteq Y$ and a mapping $g\in C(Y,X)$. If $f(g(y))=y$ on $\K$, there is nothing to prove. If not, then there is some  $ \bar{y}\in \mathcal{K}$ with $f(g(\bar{y}))\neq \bar y$. Then the mapping $u:Y\to Y$ defined by $u(y):=y-f(g(y))$   is continuous on $  Y$ and such that  $u\nequiv 0$. In the Banach space $(C(\K,X),\|\cdot\|_{\K})$ consider the function
\[
\mu(g):=\|u\|_{\K}=\sup_{y\in \K}\|u(y)\|_Y.
\]
Obviously, $\mu:C(\K,X)\to\R^+$ is a lower semicontinuous function such that $\mu (g)>0$. Roughly speaking,  the function  tells us how ''far'' the mapping   $g\in C(\K,X)$ is  from being a right inverse of $f$ on $\K$. Therefore, our task is to find a continuous mapping $\bar{g}$ from $\K$ to $X$ at which $\mu$ vanishes. We will acomplish this in several steps.

\textsc{Step 1.} For any fixed $y\in Y$ we consider the minimization problem
\[
(P) \qquad \min_{w\in X}\left\|f'(g(y))w-u(y)\right\|_Y.
\]
It is clear that a solution to $(P)$ is
\be\label{eq:w}
w(y):=[f'(g(y))]^{-1}u(y)
\ee
 at which $(P)$ attains its minimal value  zero. Furthermore, as $f'(x)$ is injective, $w(y)$ is the unique solution. In this way, one generates a single valued mapping $w:Y\to X$, defined by \eqref{eq:w}. We should note that $w$ is not continuous in general.
However, for any fixed $ y\in Y$ the function
\[
 \|f'(g(\cdot))w(y)-u(\cdot)\|_Y
\]
is continuous on $Y$ and vanishes at $y$.

Now, take  $\e>0$ such that
\be
\label{eq:EpsMu}
\e<\mu(g)/4.
\ee
For every fixed $y\in Y$, by continuity of $\|f'(g(\cdot))w(y)-u(\cdot)\|_Y$ at $y$, there exists $\delta_y>0$ such that $\fa \tilde{y}\in y+\delta_yB^\circ_Y$ we have
\be
\label{eq:cov1}
\|f'(g(\tilde{y}))w(y)-u(\tilde{y})\|_Y<\e/3.
\ee
For the compact set $\K$ we have an open cover
\[
\K\subseteq\bigcup_{y\in \K}\left(\{y\}+\delta_yB^\circ_Y\right),
\]
from which we can take a finite subcover
\[
\K\subseteq\bigcup_{i=1}^n\left(\{y_i\}+\delta_{y_i}B^\circ_Y\right)
\]
for some $\{y_1,\dots,y_n\}\subseteq \K$. Set $w_i:=w(y_i)\in X$ for $i=1,\dots, n$, and consider the finite dimensional subspace $W:=\text{span}\{w_1,\dots, w_n\}$ in $ X$.
Let $|\cdot|_W$ be a Euclidean norm on $W$. Set $L:=\max\{|w_i|^2_W,i=1,\dots,n\}>0$. Note that if $w_i=0$ for all $i$ then \eqref{eq:cov1} gives $\mu(g)=\|u\|_{\K}<\varepsilon/3$ which contradicts \eqref{eq:EpsMu}.

\medskip

Next, for any fixed $y\in Y$, we define the function $p_y:W\to\R$ by
\[
p_y(w):=\|f'(g(y))w-u(y)\|_Y+\frac{\e}{3L}|w|^2_W,
\]
and
consider the minimization problem
\[
(P_y) \qquad \min_{w\in W} p_y(w).
\]

\begin{cl}\label{cl:CoerStrictConv}
For any fixed $y\in Y$, the corresponding function $p_y:W\to\R$ is coercive and strictly convex.
\end{cl}

\begin{proof}
The function $p_y$ is coercive, as
\[
p_y(w)\ge\frac{\e}{3L}|w|^2_W\to+\infty \text{ when } |w|_W\to\infty,
\]
and $p_y$  is strictly convex because the norm  $\|\cdot\|_Y$ is convex and the squared Euclidean norm $|\cdot|^2_W$ is strictly convex.
\end{proof}

\medskip

Therefore, for any $y\in Y$,
the function $p_y:W\to\R$ attains its minimum on~$W$ at a single point which is the unique solution to the problem $(P_y)$.
Let us denote this solution by $z(y)$. Obviously, this defines a single-valued mapping $z:Y\to W$.

\begin{cl}\label{cl:z_bounded}
The set $z(\K)$ -- the image of $\K$ under $z$ -- is a bounded set in $W$, and  $p_y(z(y))<2\e/3, \ \fa y\in \K$.
\end{cl}
\begin{proof}
Since for any $w\in W$,
\[
\|f'(g(y))z(y)-u(y)\|_Y+\frac{\e}{3L}|z(y)|^2_W\le \|f'(g(y))w-u(y)\|_Y+\frac{\e}{3L}|w|^2_W,
\]
taking $w=0$ we get
\[
 \frac{\e}{3L}|z(y)|^2_W\le \| u(y)\|_Y .
\]
Since the mapping $u(\cdot)$ is continuous, it is bounded on the compact $\K$, and the latter  implies boundedness of $z$ on $\K$.

Let $y\in \K$. Then there exists $ y_i\in \K$ such that $y\in y_i+\delta_{y_i}B^\circ_Y$. Recalling that $w_i=w(y_i)\in W$, using \eqref{eq:cov1} and the definition of $L$, we get
\begin{eqnarray*}
p_y(z(y))&=&\min_{w\in W}p_y(w)=\min_{w\in W}\left\{\|f'(g(y))w-u(y)\|_Y+\frac{\e}{3L}|w|_W^2\right\}\\
 &\le&\|f'(g(y))w_i-u(y)\|_Y+\frac{\e}{3L}|w_i|_W^2<\frac{\e}{3}+\frac{\e}{3L}L=\frac{2\e}{3}.
\end{eqnarray*}
\end{proof}

\begin{cl}
\label{cl:z_cont}
The mapping $z:Y\to W$ is continuous on $Y$.
\end{cl}
\begin{proof}
 Let $\{\xi_n\}_{n=1}^\infty$ be a convergent sequence of points in $Y$ with limit point $\xi_0$. We have to show that $z(\xi_n)\to z(\xi_0)$ as $n$ tends to infinity.
 Since for any $w\in W$,
$
p_{\xi_n}(z(\xi_n))\le p_{\xi_n}(w)$,
i.e.
\begin{eqnarray}\label{zzz}
&\|f'(g(\xi_n))z(\xi_n)-u(\xi_n)\|_Y+\dfrac{\e}{3L}|z(\xi_n)|_W^2\\[8pt] \notag
&\le \|f'(g(\xi_n))w-u(\xi_n)\|_Y+\dfrac{\e}{3L}|w|_W^2, \ \fa w\in W,
\end{eqnarray}
taking $w=0$, we get
 \[
 \frac{\e}{3L}|z(\xi_n)|_W^2\le \| u(\xi_n)\|_Y
\]
 which by continuity of $u(\cdot)$ implies that $\{z(\xi_n)\}_n$ is a norm bounded sequence in the finite dimensional space $W$. Hence, we can take a convergent subsequence, $z(\xi_{n_j})\to \bar z$ as $j\to \infty$.

 Fixing $w\in W$ and using the joint continuity of $f'(\cdot,\cdot)$, and continuity of $g(\cdot)$ and $u(\cdot)$,  passing to the limit  in \eqref{zzz} as $j\to\infty$  we obtain that
\[
\|f'(g(\xi_0))\bar{z}-u(\xi_0)\|_Y+\frac{\e}{3L}|\bar{z}|^2_W\le \|f'(g(\xi_0))w-u(\xi_0)\|_Y+\frac{\e}{3L}|{w}|^2_W, \ \fa w\in W.
\]
Therefore,
\[
\|f'(g(\xi_0))\bar{z}-u(\xi_0)\|_Y+\frac{\e}{3L}|\bar{z}|^2_W\le \inf_{w\in W}\left\{\|f'(g(\xi_0))w-u(\xi_0)\|_Y+\frac{\e}{3L}|{w}|^2_W\right\}.
\]
The latter means that  $\bar{z}$ is a solution to $(P_{\xi_0})$, and since the solution is unique,  $\bar{z}\equiv z(\xi_0) $. As this is true for every choice of a convergent subsequence of $\{z(\xi_n)\}_{n=1}^\infty$, it has a single accumulation point and the continuity of $z$ follows.
\end{proof}

\medskip

\textsc{Step 2.} For $t\in \R$ we will consider the mappings $g_t:Y\to X$ defined by $g_t(y)=g(y)+tz(y)$ for $y\in Y$. By Claim \ref{cl:z_cont} it is clear that for every $t\in\R$, $g_t\in C(Y,X)$. We denote $u_t(y):=y-f(g_t(y))$ and for any fixed $y\in \K$, we consider the function $\ffi_y:\R\to \R^+$ defined by
\[
\ffi_y(t):=\|y-f(g_t(y))\|_Y.
\]
Then, for the  function $\mu$ at $g_t$ we have
\[
\mu(g_t)=\max_{y\in\K}\{\ffi_y(t)\}=\max\{\|u_t(y)\|_{Y}:\ y\in\K\}=\|u_t\|_{\K}.
\]
We apply Lemma~\ref{lem:Approx} to $f$ which by assumption is $ C^1_G$ on $X$, the compact set $K:=g(\K)$ in $  X$ (as a continuous image of a compact set), the finite dimensional space $E:=W\subseteq X$, and $r:=\sup\{\|z(y)\|_X:y\in \K\}$ (which is finite, as $z(\K)$ is bounded by Claim \ref{cl:z_bounded})    to obtain a function $\alpha(t)\to0$ as $|t|\to0$ such that for $|t|\le1$,
\begin{equation}
\label{eq:Estg_t}
\max_{y\in \K}\left\{\|f(g_t(y))-f(g(y))-tf'(g(y))z(y)\|_Y\right\}\le|t|\alpha(t).
\end{equation}

Let us  denote for simplicity $q(y):=f'(g(y))z(y)$, thus $q\in C(Y,Y)$. For any fixed $ y\in \K$, by \eqref{eq:Estg_t} we have that for $t\in[0,1]$,
\begin{eqnarray*}
\ffi_y(t)&=&\|y-f(g_t(y))\|_Y\\
         &=&
\|y-tf'(g(y))z(y)-f(g(y))-[f(g_t(y))-f(g(y))-tf'(g(y))z(y)]\|_Y\\
 &\le& \|y-tf'(g(y))z(y)-f(g(y))\|_Y+\|[f(g_t(y) )-f(g(y))-tf'(g(y))z(y)]\|_Y\\
 &\le&\|u(y)-tq(y)\|_Y+t\alpha(t).
\end{eqnarray*}
On the other hand,
\begin{eqnarray*}
\|u(y)-tq(y)\|_Y&=&\|u(y)-tu(y)+t[u(y)-q(y)]\|_Y\\
&=&\|(1-t)u(y)+t[u(y)-q(y)]\|_Y\\
 &\le& (1-t)\|u(y)\|_Y+t\|u(y)-q(y)\|_Y.
\end{eqnarray*}
Furthermore,
\begin{eqnarray}\label{u-q}
\|u(y)-q(y)\|_Y &=&\|f'(g(y))z(y)+f(g(y))-y\|_Y\notag\\
  &\le& \|f'(g(y))z(y)+f(g(y))-y\|_Y+\frac{\e}{3 L}|z(y)|^2_W\\
  &=&p_y(z(y))<2\e/3<\e, \notag
\end{eqnarray}
by Claim \ref{cl:z_bounded}. Gathering all these yields
\[
\ffi_y(t)\le(1-t)\|u(y)\|_Y+t\e +t\alpha(t),\ \fa y\in \K,\ \forall t\in[0,1].
\]
Since $\alpha(t)\to 0$, as $|t|\to 0$,  $\alpha(t)<\mu(g)/4$ for sufficiently small   $t>0$.  Then for such $t$,
\[
\max_{y\in \K}\{\ffi_y(t)\}\le (1-t)\|u\|_{\K}+t\alpha(t)+t\e,
\]
which implies
\[
\mu(g_t)\le (1-t)\mu(g) +t\alpha(t)+t\e =\mu(g)+t[ \alpha(t)+\e-\mu(g)].
\]
In addition, by the choice of $\e$, see \eqref{eq:EpsMu}, we have that $\alpha(t)+\e<\mu(g)/2$, thus $\alpha(t)+\e-\mu(g)<-\mu(g)/2$ and finally, for sufficiently small $t>0$,
\be
\label{eq:mu_est}
\mu(g_t)\le\mu(g)-t\mu(g)/2.
\ee

Since the condition \eqref{Had} holds for $f$ we have
$\|[f'(g(y))]^{-1}\|\le M$ for all $y\in Y$.
It remains to observe that, for any $y\in \K$, using \eqref{u-q},
\begin{eqnarray*}
\| z(y)\|_X &=& \left\|[f'(g(y))]^{-1}q(y)\right\|_X \\
        &\le& \|[f'(g(y))]^{-1}\|\|q(y)\|_Y\\
         &\le& M\|q(y)\|_Y \\
         &\le& M\| q(y)+u(y)-u(y)\|_Y\\
         &\le& M(\| q(y)-u(y)\|_Y+\|u(y)\|_Y)\\
         &\le&  M(\eps+\|u(y)\|_Y).
         \end{eqnarray*}
Therefore, by the choice of $\eps$, see \eqref{eq:EpsMu},
\begin{eqnarray*}
\|z\|_{\K} &=&\max_{y\in \K}\{\|z(y)\|_X\} \\
            &\le& M(\eps+\max_{y\in \K}\{\|u(y)\|_Y\})\\
            &<& M(\mu(g)/4 + \mu(g))\\
            &=& (5M/4)\mu(g).
  \end{eqnarray*}

This means that $-\mu(g)<-(4/(5M))\|z\|_{\K}$, which  substituted into \eqref{eq:mu_est} gives that for $t>0$ small  enough,
\[
\mu(g+tz)<\mu(g)-(2/(5M))\|tz\|_{\K}.
\]
We apply Lemma \ref{lem:mu} to  the Banach space $(C(\K,X), \|\cdot\|_{\K})$ and the function $\mu$ with $r=2/(5M)$  to conclude that $\mu$ vanishes at some  $\bar{g}\in C(\K,X)$, hence  $f( \bar g(y))=y$, $\fa y\in \K$. Finally, we   continuously extend $\bar{g}$ to the whole $Y$ by applying the Dugundji Extension Theorem, see e.g. \cite[Theorem 4.1]{Dugundji}.

It remains to be noted  that if $f(0)=0$ and $0\in \K$ then one can easily modify the above by considering instead of  the Banach space $C(\K,X)$, the Banach space $C_0(\K,X)$ of all continuous $g:\K\to X$ such that $g(0)=0$. It is clear that in this case $u(0)=w(0)=0$, also $z(0)=$, so $g_t(0)=0$ and all the rest works in the same way.
\end{proof}

It easily follows that a mapping satisfying the  assumptions of the proposition above, is surjective.

\begin{thm}
\label{thm:surj}
Let $(X,\|\cdot\|_X)$ and $(Y,\|\cdot\|_Y)$ be Banach spaces. Let  the mapping $f:X\to Y$ be   $ C^1_G$ on $X$,  let it  have invertible G\^ateaux derivatives   and let the Hadamard-L\'evy condition \eqref{Had} hold.
 Then $f$ is surjective.
\end{thm}

\begin{proof}
Let $y\in Y$ be arbitrary and set $\K=\{y\}$. From Proposition~\ref{prop:RightInv} it follows that there is a continuous $g\in C(Y,X)$ such that $ f(g(y))=y$, so there is $x=g(y)\in X$ such that $f(x)=y$. Therefore $f$ is surjective, i.e. $f(X)= Y$.
\end{proof}

\section{Local invertibility}\label{sec:LocInv}

We will  start by giving  an example of a $C^1_G$ mapping  $f: H \to H$ on the separable Hilbert space $H = \ell_2(\mathbb{N})$,  such that
$
  f'(0) = \mathrm{Id}$,
but $f$ is  not invertible in any neighbourhood of $0$. The example is elaborated on the basis of \cite[Example 3]{AS}.

\begin{example}\label{Example}
Let us consider the  separable Hilbert space $H = \ell_2(\mathbb{N})$. There exists a   mapping $\varphi:H\to H$ which is $C^1_G$ on $H$ such that
\begin{equation}\label{eq:phi-prime-zero}
  \varphi(0) = 0,\quad \varphi'(0) = 0,
\end{equation}
and for $  \lambda_n := 2^{-n^2}$
it holds that
\begin{equation}\label{eq:phi-fixed}
  \varphi(\lambda_n e_n) = \lambda_n e_n,
\end{equation}
where $\displaystyle
  e_n = (\underbrace{0, \dots, 1}_{n}, 0, \dots)$, $ n \in \mathbb{N}
$
is the canonical basis of $H$.

Then the mapping
$
  f(x) := x - \varphi(x)
$
is $C^1_G$ on $H$ with  $f'(0) = \mathrm{Id}$, but
$$
  f(\lambda_n e_n) = f(0) = 0, \qquad \forall n \in \mathbb{N},
$$
so $f$ is not injective in any neighbourhood of the origin.
\end{example}
\begin{proof}
Let $\psi\colon \mathbb{R} \to \mathbb{R}^+$ be a $C^\infty$-smooth bump function with
$
  \psi(0) = 1 = \max \psi$,
$
\ds  \operatorname{supp} \psi \subseteq \left[-\tfrac{1}{2}, \tfrac{1}{2}\right]$,
i.e., $\psi(t) = 0$, for all $t$ with $|t| \ge 1/2$.

For $a, b > 0$ and $h\in H$ with $\|h\|=1$ we set
$$
  \varphi_{h,a,b}(x) := \psi( \|P_h x\|/a)\, \psi(\|P_h^\perp x\|/b),
$$
where $P_h$ is the orthogonal projection onto $\mathbb{R}h$, e.g.
$P_{e_n} x = \langle x, e_n\rangle\, e_n$, and $P_h^\perp$ is the orthogonal projection onto
the hyperplane orthogonal to~$h$.

It is straightforward to check that  the mapping
$$
  \varphi(x) := \sum_{n=1}^{\infty} \lambda_n\, \varphi_{e_n,\, \lambda_n,\, 2^n\lambda_n}(x - \lambda_n e_n)\, e_n
$$
is $C^1_G$ on $H$ and satisfies \eqref{eq:phi-prime-zero} and \eqref{eq:phi-fixed}.
\end{proof}

From this example it is clear that  for a $C^1_G$ mapping one needs to assume more than just invertibility of its G\^ateaux derivative   to ensure that it is locally invertible. The following local inversion result    can be found in~\cite{IJLZ}.

\begin{thm}
\label{thm:local}\cite[Theorem 10]{IJLZ}
Let $(X,\|\cdot\|_X)$ and $(Y,\|\cdot\|_Y)$ be Banach spaces and  $f:X\to Y$  be a continuous mapping which is   G\^ateaux differentiable.

Suppose that for some $x\in X$ the derivative $f'(x) $ is invertible and
\begin{equation}\label{loc}
\limsup_{x'\to x} \| f'(x') - f'(x) \|  < \|[f'(x)]^{-1} \|^{-1}.
\end{equation}
Then there exist open sets $U \subseteq X$ and $V \subseteq Y$,  containing $x $ and $f(x )$ respectively, such that
 $ f |_{U }:U  \to V $ is a bi-Lipschitz homeomorphism, and
   $g:=[f |_{U }]^{-1}$ is G\^{a}teaux differentiable on $V $ with $g'(y)=[f'(g(y))]^{-1}$ for $y\in V $.
 If in addition $f$  is $ C^1_G$ on $U $,  then  $g $  is $ C^1_G$ on $V $.
\end{thm}

\section{Main result}\label{sec:Main}

Here we  present the proof of our main result along the lines of the proof of \cite[Theorem~1]{IZ-Hadamard}.

\begin{thm}\label{thm:Main}
Let $(X,\|\cdot\|_X)$ and $(Y,\|\cdot\|_Y)$ be Banach spaces and  $f:X\to Y$  be a $C^1_G$ mapping on $X$ such that for any $x\in X$ its G\^ateaux derivative $f'(x) $ is invertible and  \eqref{loc} holds. Let also the Hadamard-L\'evy condition \eqref{Had} hold.
Then $f$ is $C^1_G$ invertible on $X$, i.e. there exists $g:Y\to X$ which is $C^1_G$ on $Y$, such that
\[
g(f(x))=f(g(x))=x,\quad \fa x\in X.
\]
Moreover, $g$ is globally Lipschitz with   Lipschitz constant $M$.
\end{thm}

\begin{proof}
First we will show that $f$ is bijective.

As the  conditions of Theorem~\ref{thm:surj} are satisfied,   the surjectivity of $f$ follows. So, $f(X)= Y$.

To get injectivity of $f$, we exploit  the idea of the proof of \cite[Theorem 1]{IZ-Hadamard}. Let us take $a,b\in X$ be such that $f(a)=f(b)$.
By considering instead of $f$ the function $x \to f(a) - f(a-x)$
we can and do assume without loss of generality that $b=0$  and $ f(0)=0$. Then $  f(a) = 0$ and we need to show that $a=0$. Set $\K := f([0,a])$.
Since $f$ is continuous, $\K$ is a nonempty compact.

As the assumptions of Proposition~\ref{prop:RightInv} are fulfilled, there is a continuous mapping
$$
  g:Y\to X,\mbox{ such that }g(0) = 0\mbox{ and }f(g(y))=y,\quad\forall y\in \K.
$$
Consider the set
$$
  I := \{t\in [0,1]:\ g(f(ta)) = ta\}.
$$
Obviously, $0\in I$, because $g(0) = 0$. Due to the continuity of $g$ and $f$ the set $I$ is closed and, therefore, compact. Let
$$
   \bar t := \max \{ t:t\in I\}.
$$
Assume that $\bar t < 1$. By the local inversion result -- Theorem~\ref{thm:local} applied at $\bar ta$ --  there are open sets $U\ni \bar ta$ and $V\ni f(\bar ta)$ such that for each $y\in V$
there is a unique $x\in U$ satisfying
$
  f(x) = y$.

From the continuity of $f$  and $g$ there is $\beta > 0$ such that for all $t\in(\bar t, \bar t +\beta)\subseteq (0,1)$ we have $ta \in U$, $g(f(ta))\in U$, and  $f(ta) \in V$.

Hence, for   $t\in(\bar t, \bar t +\beta)$, the point $ta$ is the only solution of $f(\cdot)=f(ta)$.

Since $f(g(y))=y$ for all $ y\in \K$ and since  $f(ta)\in \K$, we have that $ f(g(f(ta))) = f(ta)$.

From the uniqueness of the solution to $f(\cdot )= f( ta)$ one gets that $g(f(ta)) = ta$ for all $t\in(\bar t, \bar t +\beta)$ which contradicts the definition of $\bar t $.

So, $\bar t =1$ meaning that $g(f(a)) = a$. But $f(a) = 0$. Since $g(0)=0$, it follows that $a=0$.

We have proved that if $f(a)=f(b)$ then $a=b$, so $f$ is injective.

Till now we have established that $f$ is invertible and for any $y\in Y$ there is a unique $x\in X$ such that $f(x)=y$. Therefore,  there is a mapping $g:=f^{-1}:Y\to X$.    From Theorem~\ref{thm:local} it follows that for any $x\in X$ there is an open set $U\ni x$ and an open set $V\ni f(x)$ such that $[f |_U]^{-1}$ is $C^1_G$ on $V$.  In other words,  $g$ is $C^1_G$ on $V$, and in particular $g$  is $C^1_G$ at $y$. As $y$ was arbitrary, $g$ is $C^1_G$ on $Y$.

\medskip

Finally,  as $g'(y)=[f'(g(y))]^{-1}$ for all $y\in Y $ the condition \eqref{Had} implies that $\|g'(y)\|\le M$ for all $y\in Y $. The global Lipschitz continuity of the $C^1_G$ mapping $g$   easily  follows from   \eqref{int2}.
\end{proof}

\bigskip

\textbf{Data Availability Statement.}
No datasets were generated or analyzed during the current study.

\bigskip

\textbf{Funding.}
M. Hamamdjiev is supported by the Science Fund of Sofia University,  grant No. 80-10-176/22.04.2026.
N. Zlateva is supported by the Bulgarian National Science Fund,  grant No. KP-06-H92/6 signed December 8, 2025.

\bigskip

\textbf{Acknowledgements.}   The authors are extremely grateful to Dr. Milen Ivanov for his constant support and encouragement and especially for the help with Example~\ref{Example}.


\begin{thebibliography}{99}

\bibitem{Aubin-Ekeland-book}
J.-P.  Aubin,  I.  Ekeland. Applied nonlinear analysis, New York, John Wiley \& Sons, 1984.

\bibitem{AubinFrankowska}
J.-P.  Aubin,  H. Frankowska, On inverse function theorems for set-valued maps, J. Math. Pures et Appl. 66 (1987)  71--89.


\bibitem{AS}
V. Averbukh, O. Smolyanov, The various definitions of the derivative in linear topological spaces,
Russian Math. Surveys 23(4)  (1968) 67--113.

\bibitem{CFR}
R. Cibulka, M. Fabian, T. Roubal, An Inverse Mapping Theorem in Fr\'echet-Montel Spaces, Set-Valued Var. Anal. 28 (2020) 195--208.


\bibitem{Dugundji}
J. Dugundji, An extension of Tietze's theorem, Pacific J. Math. 1 (1951), 353--367.

\bibitem{Ekeland}
I. Ekeland, An inverse function theorem in Fr\'chet spaces, Ann. Inst. Henri Poincar\'e (C) Anal. Non Lineaire 28(1) (2011)  91--105.

\bibitem{EkelandSere}
I.  Ekeland, E. S\'er\'e, A local surjection theorem with continuous inverse in Banach spaces, Rend. Lincei Mat. Appl. 35 (2024), 1--15.

\bibitem{Gutu}
O. Gut\'u, On global inverse theorems, Topol. Methods in Nonlinear Anal. 49(2) (2017)  401--444.

\bibitem{GutuJaramillo}
O. Gut\'u, J. A. Jaramillo, Surjection and inversion for locally Lipschitz maps between Banach spaces, J. Math. Anal. Appl. 478(2) (2019) 578--594.

\bibitem{Hadamard}
J. Hadamard, Sur les transformations ponctuelles, Bull. Soc. Math. France 34 (1906)  71--84.

\bibitem{Hamilton}
R. S. Hamilton,  The inverse function theorem of Nash and Moser,  Bull. Amer. Math. Soc. (N.S.)  7(1) (1982)  65--222.


\bibitem{Ioffe-book}
A.  Ioffe, Variational Analysis of Regular Mappings: Theory and Applications, Springer Monogr. Math. Cham, Springer, 2017.

\bibitem{IZ-Hadamard}
M. Ivanov,  N. Zlateva. Hadamard inverse function theorem proved by variational analysis, Serdica Math. J. 49(1-3) (2023)  97--106.

\bibitem{IJLZ}
M. Ivanov, J. A. Jaramillo,  S. Lajara,  N. Zlateva, Continuous selections and invertibility of nonsmooth maps between Banach spaces, J. Optim. Theory Appl.   206 (2025), art. no. 41.

\bibitem{Levy}
P. L\'evy, Sur les fonctions des lignes implicites, Bull. Soc. Math. France 48 (1920) 13--27.

\bibitem{Pales}
Z. P\'ales, Inverse and Implicit Function Theorems for Nonsmooth
Maps in Banach Spaces, J. Math. Anal. Appl. 209(1) (1997) 202--220.

\bibitem{Plastock}
R. Plastock, Homeomorphisms between Banach spaces, Trans. Amer. Math. Soc.  200 (1974) 169--183.

\bibitem{Pourciau1}
B. H. Pourciau, Hadamard's theorem for locally Lipschitzian maps, J. Math. Anal. Appl. 85(1) (1982) 279--285.

\bibitem{Pourciau2}
B. H. Pourciau, Global invertibility of nonsmooth mappings, J. Math. Anal. Appl. 131(1) (1988) 170--179.


\bibitem{Takahashi}
W. Takahashi,   Existence theorems generalizing fixed point theorems for multivalued mappings, In: M. A. Thera, J. B. Baillon (Eds.)  Fixed point theory and applications, Pitman Research Notes in Mathematics Series   252 (1991), Longman Scientific \& Technical, Harlow,   397--406.



\end{thebibliography}
\end{document}